\documentclass[11pt]{amsart}

\usepackage{epsfig, amsmath, amssymb, amsthm, amsfonts}

\newcommand{\df}{\ensuremath{\partial}}
\newcommand{\alg}{\ensuremath{\mathcal{A}}}
\newcommand{\rr}{\ensuremath{\mathbb{R}}}
\newcommand{\zz}{\ensuremath{\mathbb{Z}}}

\theoremstyle{plain}
\newtheorem{thm}{Theorem}
\newtheorem{cor}[thm]{Corollary}
\newtheorem{lem}[thm]{Lemma}

\theoremstyle{definition}
\newtheorem{defn}[thm]{Definition}

\newcommand{\Aug}{\operatorname{Aug}}
\newcommand{\sgn}{\operatorname{sgn}}

\begin{document}

\title[Augmentations and Rulings for Legendrian Knots]{The
Correspondence between Augmentations and Rulings for Legendrian
Knots}

\author[L. Ng]{Lenhard L. Ng}
\address{Stanford University, Stanford, CA 94305}
\thanks{LLN is supported by the American Institute of Mathematics.}

\author[J. Sabloff]{Joshua M. Sabloff}
\address{Haverford College, Haverford, PA 19041}

\begin{abstract}
  We strengthen the link between holomorphic and generating-function
  invariants of Legendrian knots by establishing a formula relating
  the number of augmentations of a knot's contact homology to
  the complete ruling invariant of Chekanov and Pushkar.
\end{abstract}

\maketitle

\section{Introduction}
\label{sec:intro}

The theory of Legendrian knots in the standard contact $\rr^3$ up to
Legendrian isotopy is a subject of much recent interest. Beginning
in the late 1990's, two ``non-classical'' invariants of such
Legendrian knots emerged. The first is based on the contact homology
theory of Eliashberg and Hofer \cite{bib:El}, and was couched in
combinatorial form by Chekanov \cite{bib:Ch}.  This invariant is a
differential graded algebra that counts certain holomorphic disks in
the symplectization of $\rr^3$. The second is the ruling invariant
of Chekanov and Pushkar \cite{bib:CP}, which is derived from the
theory of generating functions but has a simple combinatorial
definition. See \cite{bib:Ch2} for a brief introduction to both
invariants.

It seems likely that the information contained in the
Chekanov--Pushkar ruling invariant is in some way the same as that
contained in the linear level of the contact homology DGA. The first
link between the two theories was provided by Fuchs \cite{bib:Fu},
who introduced the concept of a ruling independently of
Chekanov--Pushkar and demonstrated that the existence of a ruling
implies the existence of a so-called augmentation of the contact
homology DGA. Fuchs and Ishkhanov \cite{bib:FI}, and, independently,
the second author \cite{bib:S} proved the converse. In fact, a more
general sort of correspondence seems to hold between
generating-function invariants and linearized contact homology; see
\cite{bib:NT} for evidence in $J^1(S^1)$.

In this paper, we deepen the link between contact homology and rulings
by showing that there is a many-to-one correspondence between
augmentations and rulings for a plat-position front of a Legendrian
knot.\footnote{This is in fact true for any Legendrian front, but we
  restrict to plat position for ease of exposition.} An algorithm was
given in \cite{bib:S} for obtaining a ruling of a plat-position front
from an augmentation of its DGA.  We show that the number of
augmentations which yield a particular ruling via the algorithm is
determined by the combinatorics of the ruling (Theorem~\ref{thm}).
This shows, in particular, that the total number of augmentations of
the DGA of any Legendrian knot can be calculated from its ruling
invariant (Corollary~\ref{cor}).

\section{Main Results}
\label{sec:results}

We begin by introducing the terminology necessary to state our main
results. Please refer to \cite{bib:S} (or any number of other papers)
for the precise definitions and conventions which we will use: contact
homology DGA, stable tame isomorphism, and algebraic stabilization. As
some geometric conventions and the augmentation and ruling
constructions are central to this paper, we will review them in more
detail.

\begin{figure}
  \centerline{\includegraphics{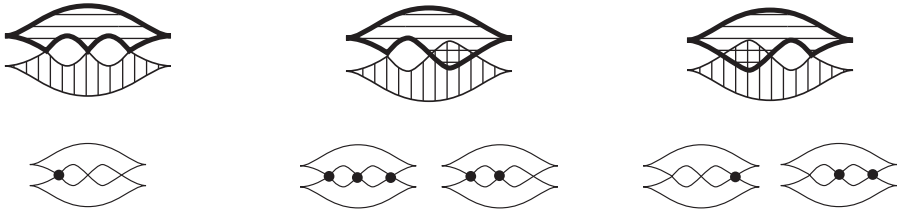}}
  \caption{The many-to-one correspondence between $\zz$-graded augmentations
    and rulings on a Legendrian trefoil knot. Large dots represent augmented
    crossings. The trefoil pictured has $r(K)=0$ and $\chi^*_0(\alg)
    = 1$; the leftmost ruling has $\theta=-1$, while the other two
    have $\theta=1$.
    }
  \label{fig:trefoil-ex}
\end{figure}

A \textbf{front diagram} of a Legendrian knot $K$ in the standard
contact space $(\rr^3, dz - y\, dx)$ is simply the image of the
projection of $K$ to the $xz$ plane.  A front diagram is in
\textbf{plat position} if all of the left cusps have the same $x$
coordinate, all of the right cusps have the same $x$ coordinate, and
no two crossings have the same $x$ coordinate.  For example, the front
diagram of the Legendrian trefoil in Figure~\ref{fig:trefoil-ex} is in
plat position.

Next, we describe some algebraic constructs related to the contact
homology DGA. Let $K$ be a Legendrian knot in the standard contact
$\rr^3$.  Its contact homology is the homology of a semifree unital
differential graded algebra $(\alg, \df)$ with coefficients in
$\zz/2$ and grading over $\zz/(2r(K))$.  Here, $r(K)$ is the
rotation number of $K$, and $\alg$ is generated by the crossings and
right cusps of a front diagram $D$ of $K$. Of particular interest
for this paper are the gradings of the crossings: define a locally
constant function $\mu$ from the complement of the cusps in $K$ to
$\zz/(2r(K))$ such that at each cusp, $\mu$ increases by $1$ from
the lower to the upper strand.  Near a crossing $q$, let $\alpha$
(resp.\ $\beta$) be the strand of $D$ with more negative (resp.\
positive) slope.  Assign the grading $|q| \equiv \mu(\alpha) -
\mu(\beta) \pmod{2r(K)}$.  For instance, all three crossings in the
trefoil of Figure~\ref{fig:trefoil-ex} have degree $0$. (The other
generators of $\alg$, the right cusps, are all assigned grading
$1$.) The central result of \cite{bib:Ch} states that $(\alg,\df)$
is invariant under Legendrian isotopy, up to stable tame
isomorphism.

Henceforth we assume that we have a integer $\rho$ with
$\rho\,|\,2r(K)$. In this case, it makes sense to discuss whether a
crossing of $D$ (a generator of $\alg$) has degree divisible by
$\rho$; we write $|a|\in\zz/(2r(K))$ for the degree of $a$.

One way to extract information from $(\alg,\df)$ is via
augmentations. A \textbf{$\pmb{\rho}$-graded augmentation} of
$(\alg,\df)$ is an algebra map $\varepsilon:\thinspace\alg\to\zz/2$
such that $\varepsilon\circ\df = 0$, $\varepsilon(1)=1$, and
$\varepsilon(a)=0$ if $\rho \nmid |a|$. (This last condition states
that $\varepsilon$ respects the grading of $\alg$ over $\zz/\rho$,
where the base ring $\zz/2$ is considered to lie in degree $0$.) Two
cases are of special interest. If $\rho=1$, then a $1$-graded
augmentation is also known as an \textbf{ungraded augmentation}, on
which any generator can be augmented (i.e., sent to 1 by
$\varepsilon$) regardless of grading. If $\rho=0$, and hence
$r(K)=0$, then a $0$-graded augmentation is also known as a
\textbf{$\zz$-graded augmentation}, on which only generators of
degree $0\in\zz$ can be augmented.

Note that, since $\alg$ is finitely generated as an algebra, any DGA
has only finitely many $\rho$-graded augmentations. The number of
$\rho$-graded augmentations of $(\alg,\df)$ is invariant under tame
isomorphism, and changes by a power of $2$ under algebraic
stabilization. More precisely, an algebraic stabilization of degree
$i$ adds generators $\alpha,\beta$ of degree $i-1,i$ respectively,
with $\df(\beta)=\alpha$ and $\df(\alpha)=0$. It follows that this
stabilization leaves the number of $\rho$-graded augmentations
unchanged if $\rho\nmid i$, and doubles this number if $\rho\,|\,i$.

If $\rho$ is odd or $\rho=0$, a normalized count of the number of
$\rho$-graded augmentations gives us an invariant of $(\alg,\df)$
under stable tame isomorphism. Let $a_k$ be the number of generators
of \alg\ of degree $k \pmod{\rho}$; the set $\{a_k\}_{k \in
\zz/\rho}$ is called the \textbf{degree distribution} of the DGA.
Set $\chi^*_\rho(\alg)$ to be the following ``shifted Euler
characteristic'' of \alg:
\[
  \chi^*_0(\alg) = \sum_{k\geq 0} (-1)^k a_k + \sum_{k<0}
  (-1)^{k+1} a_k
\]
and
\[
  \chi^*_\rho(\alg) = \sum_{k=0}^{\rho-1} (-1)^k a_k
\]
if $\rho$ is odd.  We can now define the following invariant of
Legendrian isotopy:

\begin{defn}
  Let $K$ be a Legendrian knot and $\rho$ a divisor of $2r(K)$ with
  $\rho=0$ or $\rho$ odd. The \textbf{normalized $\pmb{\rho}$-graded
  augmentation number} of $K$, written $\Aug_\rho(K)$, is the
  product of $2^{-\chi^*_\rho(\alg)/2}$ and the number of
  $\rho$-graded augmentations of $(\alg,\df)$, where $(\alg,\df)$ is
  the contact homology DGA of $K$.
\end{defn}

We remark that $\Aug_0$ distinguishes between Chekanov's well-known
examples of Legendrian $5_2$ knots. Note that, when $\rho$ is even
and nonzero, it is not possible to define similar normalized
augmentation numbers to give an invariant under stable tame
isomorphism. The issue in this case is that one can add an algebraic
stabilization in each of degrees $0,2,\dots,\rho-2$, or in each of
degrees $1,3,\dots,\rho-1$, and the degree distribution changes the
same way in both cases; however, the number of augmentations doubles
in the former case but remains the same in the latter.

We next turn to rulings. A \textbf{ruling}\footnote{Chekanov
  \cite{bib:Ch2} calls this an ``admissible decomposition'' and
  Chekanov and Pushkar \cite{bib:CP} call this a ``positive proper
  decomposition''.} of a front is a decomposition of the front into
pairs of paths beginning at a left cusp and ending at a right cusp
that co-bound disks. Each ruling path is smooth except at cusps and
certain crossings called switches. Near a switch, the ruling disks
must look like a diagram from the first row of
Figure~\ref{fig:crossing-config}. If $\rho\,|\,2r(K)$, then the
ruling is \textbf{$\pmb{\rho}$-graded} if all switches have degree
divisible by $\rho$.

To each ruling $R$ of a front $D$ of a Legendrian knot $K$, we can
associate an integer $\theta(R) = c(D) - s(R)$, where $c(D)$ is the
number of right cusps of $D$ (i.e., the number of closed curves in
$R$), and $s(R)$ is the number of switches of $R$. A main result of
Chekanov and Pushkar \cite{bib:CP} is that the multiset
\[
  \Theta_\rho(K) = \{\theta(R)\,|\, R \text{ is a $\rho$-graded ruling
  of } D\},
\]
which we call the \textbf{complete ruling invariant} of $K$, is
invariant under Legendrian isotopy.

The complete ruling invariant is effective: Chekanov \cite{bib:Ch2}
used the complete ruling invariant to distinguish his original $5_2$
examples of Legendrian knots that have the same classical
invariants. We conjecture in passing that when $\rho=1$ (i.e., the
ungraded case), it seems possible that $\Theta_1(K)$ always depends
only on the topological type and Thurston--Bennequin number of $K$.
If this were true, then $\Theta_1$ would provide a topological knot
invariant by considering any Legendrian $K$ in the knot type with
maximal Thurston--Bennequin number. This is related to, but
independent of, the conjecture in \cite{bib:Ng} that the ungraded
abelianized characteristic algebra is a topological invariant, and
the ``irresponsible conjecture'' in \cite{bib:Fu} (see also
\cite{bib:FI}) that $\Theta_1(K) = \emptyset$ if and only if the
Kauffman polynomial bound on the Thurston--Bennequin number of $K$
is not sharp.

We are now in a position to state our main results.

\begin{thm} \label{thm}
  Given a Legendrian knot $K$ with a diagram $D$ in plat position, and
  an integer $\rho$ with $\rho\,|\,2r(K)$ and either $\rho=0$ or
  $\rho$ odd, there is a many-to-one correspondence between
  $\rho$-graded augmentations of the contact homology DGA $(\alg,
  \df)$ of $D$ and $\rho$-graded rulings of $D$, with
  $2^{(\theta(R) + \chi^*_\rho(\alg))/2}$ augmentations
  corresponding to a ruling $R$.
\end{thm}

\noindent An example of the correspondence appears in
Figure~\ref{fig:trefoil-ex}.

Theorem~\ref{thm} has the following immediate consequence.

\begin{cor} \label{cor}
  Suppose $\rho \,|\,2r(K)$ and $\rho$ is zero or odd.  The
  normalized $\rho$-graded augmentation number of a Legendrian knot
  $K$ can be deduced from the complete ruling invariant:
  \begin{equation*}
    \Aug_\rho(K) = \sum_{\theta \in \Theta_\rho(K)}
    2^{\theta/2}.
  \end{equation*}
\end{cor}

As stated earlier, it is possible that the complete ruling invariant
is actually contained somehow in the contact homology DGA. One could
hope to associate to each augmentation a fractional power of $2$, so
that if a ruling gives $2^k$ augmentations then each of those
augmentations is given the fraction $2^{-k}$. This would allow us to
recover the complete ruling invariant from the augmentations of the
DGA.

\section{Proofs}
\label{sec:proofs}

The proof of Theorem~\ref{thm} rests on an algorithm described in
\cite{bib:S} that produces a ruling from an augmentation. After
setting notation, we will describe this algorithm, use it to define
the many-to-one correspondence in the theorem, and finally count the
number of augmentations that correspond to each ruling in terms of
the degree distribution of the contact homology DGA.

\begin{figure}
  \centerline{\includegraphics[height=3.5in]{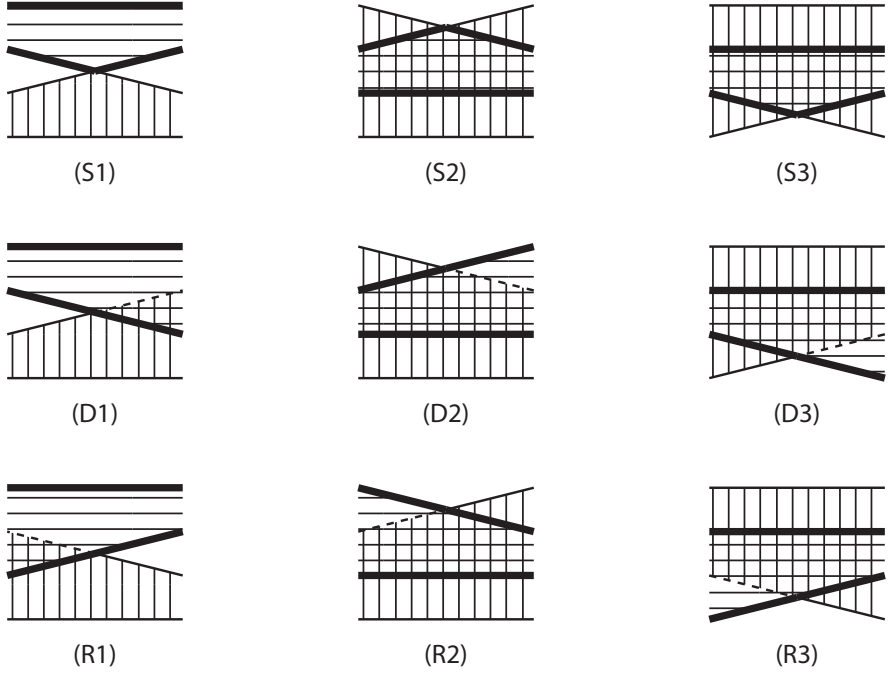}}
  \caption{The configurations labeled (S*), (D*), and (R*) are
    switches, departures, and returns, respectively.}
  \label{fig:crossing-config}
\end{figure}

Let $K$ be a Legendrian knot with plat diagram $D$ and a $\rho$-graded
ruling $R$. Near a crossing, call the two ruling paths that are
incident to the crossing \textbf{crossing paths} and call the ruling
paths that are paired with the crossing paths \textbf{companion
  paths}.  At a crossing of $D$ whose grading is not divisible by
$\rho$, the ruling simply passes through the crossing with no switch.
We can use the ruling $R$ to partition the crossings of $D$ with
grading divisible by $\rho$ into three types:
\begin{description}
\item[Switches] The ruling disks are nested or disjoint on both sides
  of the crossing. See the first line of
  Figure~\ref{fig:crossing-config}.
\item[Departures] Moving left-to-right, the ruling disks pass from
  nested or disjoint to interlaced.  See the second line of
  Figure~\ref{fig:crossing-config}.
\item[Returns] Moving left-to-right, the ruling disks pass from
  interlaced to nested or disjoint.  See the third line of
  Figure~\ref{fig:crossing-config}.
\end{description}

\noindent Henceforth the terms ``switches'', ``departures'', and
``returns'' will refer to crossings of grading divisible by $\rho$
with the above properties.

The correspondence in Theorem~\ref{thm} is derived from the procedure
implicitly described in \cite{bib:S}\footnote{See \cite[\S
  3.3]{bib:S}, especially the remark at the end of the section.} for
producing a ruling on a plat diagram from an augmentation
$\varepsilon$.  We recall the procedure here. Label the crossings
from left to right by $q_1, \ldots, q_n$. The construction begins at
the left cusps, where any ruling must pair paths incident to the
same cusp. The next step is to extend the ruling over the crossings
from left to right, one crossing at a time.  The extension procedure
uses a succession of ``virtual augmentations'', $\rho$-graded
algebra maps $\varepsilon_j:\thinspace\alg\to\zz/2$, where
$\varepsilon_1 = \varepsilon$ and $\varepsilon_{j+1}$ is defined
after the ruling has been extended over $q_j$.

The extension of the ruling over $q_j$ depends on the configuration
of the ruling disks to the left of $q_j$ and on
$\varepsilon_j(q_j)$. If $|q_j|$ is not divisible by $\rho$, then
the ruling extends without switching.  Otherwise, the ruling extends
according to the following instructions:
\begin{itemize}
\item If the ruling disks are nested or disjoint to the left of $q_j$
  and $\varepsilon_j(q_j) = 1$, then extend the ruling over the
  crossing as a switch.
\item If the ruling disks are nested or disjoint to the left of $q_j$
  and $\varepsilon_j(q_j) = 0$, then extend the ruling over the
  crossing as a departure.
\item Otherwise --- i.e., if the ruling disks are interlaced to the
  left of $q_j$ --- then extend the ruling over the crossing as a
  return.
\end{itemize}

Once the ruling has been extended over $q_j$, we define the virtual
augmentation $\varepsilon_{j+1}$. If $\varepsilon_j(q_j)=0$, or if the
ruling is extended over $q_j$ as in configuration (R1), then set
$\varepsilon_{j+1}=\varepsilon_j$. For the other configurations, set
$\varepsilon_{j+1}(q_k)=\varepsilon_j(q_k)$ for $k\leq j$; we now
inductively define $\varepsilon_{j+1}(q_k)$ for $k\geq j+1$.

\begin{figure}
  \centerline{\includegraphics[height=0.8in]{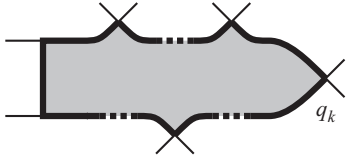}}
  \caption{The embedding of the disk used to determine the new virtual
    augmentation $\varepsilon_{j+1}$ on $q_k$.  The vertical segment
    lies to the right of $q_j$ and joins either the crossing strands or
    the companion strands; in the latter case, some other strands may
    cross the vertical segment.  The corners on the top and bottom lie at
    crossings where $\varepsilon_{j+1}=1$.}
  \label{fig:special-disk}
\end{figure}

Suppose that $\varepsilon_{j+1}(q_l)$ has been defined for $l<k$.
The value of $\varepsilon_{j+1}(q_k)$ depends on a count of a
special set of embeddings of the $2$-disk into the diagram, up to
smooth reparametrization. If the number of these disks is odd, then
$\varepsilon_{j+1}(q_k) = \varepsilon_j(q_k)+1$; the two virtual
augmentations agree on $q_k$ otherwise.  The image of the boundary
of the disks leaves $q_k$ on the upper left strand, travel leftwards
(possibly with convex corners at vertices with
$\varepsilon_{j+1}=1$) to the top of a vertical line segment that
joins two strands of the diagram and lies just to the right of
$q_j$.  The boundary then traverses this segment and travels
rightwards (again with possible convex corners where
$\varepsilon_{j+1}=1$) back to $q_k$. See
Figure~\ref{fig:special-disk}. In the case where the ruling has been
extended over $q_j$ in configurations (R2) or (R3), the vertical
segment above must join the companion strands. In the case of
configuration (S1), the vertical segment joins the crossing strands.
In the case of configuration (S2), first perform the entire
inductive procedure using disks with vertical segments that join the
companion strands, to obtain a virtual augmentation
$\varepsilon_{j+1}'$ from $\varepsilon_j$; then repeat the procedure
with vertical segments that join the crossing strands, to obtain the
virtual augmentation $\varepsilon_{j+1}$ from $\varepsilon_{j+1}'$.
For configuration (S3), do the same thing in the opposite order. It
was proven in \cite{bib:S} that this procedure determines a genuine
ruling.

With the algorithm now in hand, we may begin the proof of the
correspondence between augmentations and rulings.  The algorithm sends
an augmentation to a ruling; the question is which augmentations get
sent to the same ruling.  The answer is given by the following:

\begin{lem}
  \label{lem1}
  If a $\rho$-graded ruling of a plat diagram $D$ has $r$ returns,
  then the number of $\rho$-graded augmentations corresponding to it
  is $2^r$ if $\rho \neq 1$, and $2^{r+c(D)}$ if $\rho=1$.
\end{lem}

\begin{proof}
  The algorithm given above for obtaining a ruling $R$ from an
  augmentation $\varepsilon$ also gives a ``final'' virtual
  augmentation $\varepsilon_n$. Different augmentations result in
  different final virtual augmentations; simply consider the leftmost
  crossing $q_k$ at which the augmentations differ. It is also clear
  from the algorithm that $\varepsilon_n$ augments the switches of
  $R$, some subset of the returns of $R$, and no other crossings.

  We claim that any virtual augmentation which augments only the
  switches and some subset of the returns of $R$, for some ruling $R$,
  is the final virtual augmentation for some augmentation of the
  diagram.  Indeed, the argument in \cite{bib:S} leading to the
  algorithm shows that such a virtual augmentation $\varepsilon'$
  corresponds to an honest augmentation on the DGA $(\alg',\df')$ for
  a different (``dipped'') diagram for the knot. The algorithm
  determines a stable tame isomorphism between the DGA $(\alg,\df)$
  for $D$ and $(\alg', \df')$.  Thus, the augmentation on
  $(\alg',\df')$ yields an augmentation $\varepsilon$ on a
  stabilization of $(\alg,\df)$ and hence, by restriction, on
  $(\alg,\df)$ itself.  Then $\varepsilon'$ is the final virtual
  augmentation for $\varepsilon$, proving the claim.

  Now given a $\rho$-graded ruling $R$, the number of virtual augmentations
  augmenting only its switches and some of its returns is clearly
  $2^r$ when $\rho\neq 1$, and $2^{r+c(D)}$ when $\rho=1$, since in
  the latter case any subset of the right cusps, which have degree
  $1$, can be augmented. By the preceding argument, this is also the
  number of augmentations corresponding to $R$.
\end{proof}

Theorem~\ref{thm} now follows from Lemma~\ref{lem1} and the
following combinatorial result.

\begin{lem}
\label{lem2}
The number of returns of a $\rho$-graded ruling $R$ is
$\frac{1}{2}(\theta(R)+\chi^*_\rho(A))$ if $\rho=0$ or $\rho\geq 3$ is
odd, and $\frac{1}{2}(\theta(R)+\chi^*_1(A))-c(D)$ if $\rho=1$.
\end{lem}

\begin{proof}
  As before, let $r$ be the number of returns of a $\rho$-graded
  ruling $R$; also, let $s$ and $d$ denote the number of switches and
  departures, respectively (by definition, this only counts crossings
  of degree divisible by $\rho$).

  Consider a vertical line which intersects the plat diagram
  generically (i.e., it does not pass through a crossing or cusp).
  Since there are $c(D)$ right cusps, any such line intersects the
  diagram in $2c(D)$ points, ordered from top to bottom, and the
  ruling determines a pairing of these points. In addition, each point
  $p$ carries a Maslov index $m(p)$ given by the Maslov index of its
  strand (see \cite{bib:Ch2} or \cite{bib:S} for the definition of
  Maslov index), and the Maslov indices of paired points differ by
  $1$. Say that two pairs of points are \textbf{interlaced} if we
  encounter the pairs alternately as we move from top to bottom; that
  is, they appear from top to bottom as $a_1b_1a_2b_2$, where $a_i$
  denotes one pair of companion strands and $b_i$ denotes the other.
  To any generic vertical line we associate a number called the
  \textbf{interlacing number}. We will see that as we move the
  vertical line from left to right, beginning just after the left
  cusps and ending just before the right cusps, then the interlacing
  number begins and ends at $0$, and changes by $\pm 1$ when passing
  through each crossing; counting these changes yields the lemma.

  We first consider the case $\rho=1$, i.e., the ungraded case. Here
  we define the interlacing number of a vertical line to be the number
  of interlaced pairs on it. As we sweep the vertical line from left
  to right, the interlacing number clearly begins and ends at $0$, is
  unchanged at each switch, increases by $1$ at each departure, and
  decreases by $1$ at each return. Thus $d=r$; the lemma follows from
  the fact that $\chi^*_1(A) = a_0 = s+d+r+c(D)$.

  Next consider the case $\rho=0$, i.e., the $\zz$-graded case. For
  $k\in\zz$, let $\tilde{a}_k$ be the number of crossings of $D$ of
  degree $k$, and note that, for $a_k$ defined as in
  Section~\ref{sec:intro}, we have $\tilde{a}_k = a_k$ for $k\neq 1$
  and $a_1-\tilde{a}_1 = c(D)$; also, $a_0 = s+d+r$ and $\theta(R) =
  c(D) - s$. The lemma now reduces to proving that
  \begin{equation}
  \tag{$*$}
    d-r+\sum_{k>0}(-1)^k \tilde{a}_k + \sum_{k<0} (-1)^{k+1} \tilde{a}_k
    = 0.
    \label{eq1}
  \end{equation}

  To two interlaced pairs $a_1b_1a_2b_2$ we associate a sign:
  \[
  \sgn(a_1b_1a_2b_2) = \begin{cases} 1 & \text{if } m(a_2) - m(b_2)
    \in
    \{\ldots,-4,-2,0,1,3,5,\ldots\} \\
    -1 & \text{if } m(a_2) - m(b_2) \in
    \{\ldots,-5,-3,-1,2,4,6,\ldots\}.
  \end{cases}
  \]
  Define the interlacing number of a vertical line to be the sum over
  all interlaced pairs of this sign. Again, if we move the vertical
  line from left to right, the interlacing number begins and ends at
  $0$, and changes only when the line passes through a non-switch
  crossing.

  Consider any crossing $C$ besides a switch. We claim that the
  interlacing number changes at $C$ by $\pm 1$, and that each crossing
  contributes the appropriate sign to obtain (\ref{eq1}). If we ignore
  degree, $C$ looks like one of (D*) or (R*) from
  Figure~\ref{fig:crossing-config}. There is a one-to-one
  correspondence between interlaced pairs as for switches, except for
  one extra interlaced pair on the right of $C$, in the case of (D*),
  or on the left, in the case of (R*). Label this interlaced pair by
  $a_1b_1a_2b_2$. As we pass through $C$, the interlacing number
  changes by $\sgn(a_1b_1a_2b_2)$ for (D*) and $-\sgn(a_1b_1a_2b_2)$
  for (R*).

  Recall from \cite{bib:Ch2} or \cite{bib:S} that the degree $|\cdot|$
  of a crossing is the difference of the Maslov indices of the
  crossing strands, more precisely $m(a)-m(b)$ where $a$ lies above
  $b$ to the left of the crossing.  Note also that for the interlaced
  pair $a_1b_1a_2b_2$ as above, we have $m(a_1)=m(a_2)+1$ and
  $m(b_1)=m(b_2)+1$. It follows that the degree of a crossing that
  looks like (D1) is $m(a_2)-m(b_2)-1$; (D2) and (D3),
  $m(b_2)-m(a_2)$; (R1), $m(b_2)-m(a_2)+1$; (R2) and (R3),
  $m(a_2)-m(b_2)$. An easy computation using the definition of
  $\sgn(a_1b_1a_2b_2)$ shows that the interlacing number changes by:
  \begin{itemize}
  \item $+1$ if $|C| \in \{\ldots,-5,-3,-1,2,4,6,\ldots\}$, or $|C| =
    0$ and $C$ is a departure,
  \item $-1$ if $|C| \in \{\ldots,-6,-4,-2,1,3,5,\ldots\}$, or $|C| =
    0$ and $C$ is a return.
  \end{itemize}
  Since the interlacing begins on the left as $0$ and ends on the
  right as $0$, counting each of these changes of $\pm 1$ yields
  (\ref{eq1}), as desired. This completes the proof when $\rho=0$.

  The proof of the lemma when $\rho\geq 3$ is odd is identical to the
  proof when $\rho=0$, except that we define
  \[
  \sgn(a_1b_1a_2b_2) = \begin{cases} 1 & \text{if } m(a_2) - m(b_2)
    \in
    \{0,1,3,5,\ldots,\rho-2\} \\
    -1 & \text{if } m(a_2) - m(b_2) \in \{2,4,6,\ldots,\rho-1\}.
  \end{cases}
  \qedhere
  \]
\end{proof}


\end{document}